# Chance constrained optimization of energy intensive production as beneficial power units

*Johannes Nicklaus[1], Lea Brass[1], Gunnar Schubert[1]*

[1]HTWG Konstanz, Konstanz, Germany
*{jnicklau, lbrass, gschuber}@htwg-konstanz.de

**Keywords:** ENERGY-INTENSIVE INDUSTRY, STOCHASTIC PROGRAMMING, UNCERTAINTY ANALYSIS, CHANCE-CONSTRAINED OPTIMIZATION,

## Abstract

We study linear policy approximations for the risk-conscious operation of an industrial energy system with uncertain wind power, significant and variable electricity demand, and high thermal output, as found in a modern foundry. The system incorporates thermal storage and operates under rolling forecasts, leading to a sequential decision-making framework. To address uncertainty in key parameters, we formulate chance-constrained optimization problems that limit the probability of critical constraint violations—such as unmet demand requirements or the exceedance of system boundaries. To reduce computational effort, we replace direct uncertainty handling with a parameter-modified cost function that approximates the underlying risk structure. We validate our method through a numerical case study, demonstrating the trade-offs between operational efficiency and reliability in a stochastic environment.

## 1 Introduction

To mitigate global warming, it is imperative to substantially reduce greenhouse gas emissions. In sectors characterized by high energy demand, awareness of the associated environmental impacts is increasing, and various pathways for transitioning away from carbon-based fuels—such as the adoption of biofuels, hydrogen, and electricity—are being actively explored. Nevertheless, the integration of renewable electricity remains challenging due to significant uncertainties in its generation. These fluctuations clash with the inflexible and continuous power demands as well as the rigid production schedules prevalent in these sectors.

An effective strategy to address part of this issue involves utilizing waste heat from industrial production facilities as a thermal energy source for local district heating networks. This approach is particularly advantageous, as the periods of peak heat demand often coincide with times when renewable energy generation, particularly from photovoltaic and wind sources, is low, such as during cold winter days. Although both heat demand and renewable energy generation are influenced by some similar meteorological phenomena, their specific temporal patterns and sources of uncertainty differ significantly. Consequently, the task of aligning production schedules with fluctuating energy availability while ensuring reliable, year-round operation constitutes a complex optimization problem. To ease the task of matching production heat output with heat demand, thermal energy storage (TES) facilities can be built up. This, together with electrical energy storage (EES) for storing high peaks of renewable power generation for later use, allows for a more flexible and economical operation. However, this requires the optimization of the power flow along with the utilization of the available storage capacities.

In particular, few models capture the dynamic interaction between decentralized production facilities, fluctuating energy availability, and flexible heat markets within a unified optimization framework.

In this paper, we take a step towards addressing these challenges associated with planning the energy flow in such networks. We build upon results and methods by [1] as well as [2], by simulating renewable electricity production (photovoltaic and wind) based on realistic conditions for the German energy market. We model the energy demand of an industrial iron foundry and the energy demand of a surrounding urban area. To reflect economic incentives for heat and electricity provision, a dynamic market price for heat and electricity demand is incorporated. To optimize production and energy flows under uncertainty, we formulate a stochastic optimization problem and solve it using chance-constrained linear programming. The optimal size of EES and TES is investigated to assess the benefits of storage flexibility on overall system performance.

### 1.1 Related Literature

In recent years a lot of scientific effort has been put into optimizing energy networks in various aspects. Besides [1] and [2], optimization has been performed for the operation of power networks. For a comprehensive review of optimization methods under uncertainty [3].

In the literature deterministic approaches are prevalent and have been used to optimize the production scheduling of combined heat and power plants via scenario optimization [4], jointly optimize the production schedule of a heat and power network through distributionally robust optimization [5] or multi-objective optimization [6].

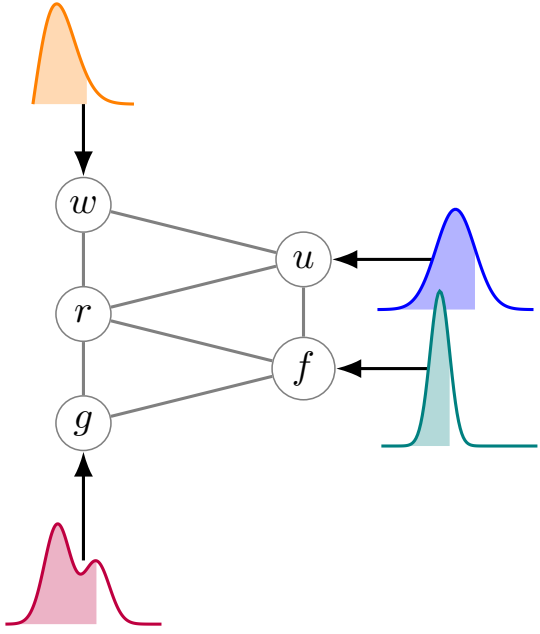

Fig. 1 For certain values influencing the energy network, we calculate probability distributions with given PDF/CDF. We then use either the mean (not shown) or a confidence level (PDF fill level above) to acquire deterministic values with which we perform the LP-optimization. In this paper we only evaluate Gaussian normal distributions.

However, uncertainty aware probabilistic approaches are still rather underexplored. In [7], the authors specifically investigated the uncertainty in heat demand and how to use it in seasonal storage. Methods of designing and using thermal heat storages in a socio-ecological and economical way have been explored in [8] and [9].
In recent years probabilistic forecasting community has been focusing increasingly on methods for electricity load forecasting [10] and renewable power generation [11] as well as the effects on the system uncertainty due to increasing penetration of renewable energy sources [12]. Dealing with this uncertainty often leads to methods from control [13], operations research [14] or optimization theory [15].
Energy systems with a large share of distributed energy storage and conversion facilities in highly connected heat and power networks have only recently drawn attention [16–18].

*1.2 Organization of this paper*

The remainder of this paper is organized as follows. section 2 introduces the notation and the theory of chance-constrained optimization. Here, the energy model is also explained. subsection 2.4 describes how uncertainty affects the decision-making process by influencing policies. subsection 2.5 explains the details of the simulation setup and the models that were used to obtain probabilistic forecasts. section 3 presents the experimental results and gives design recommendations for the investigated system. section 4 gives concluding remarks and puts the presented research in perspective.

## 2 Theory

*2.1 Notation*

The notations in this paper are fairly standard for the optimization community. All vectors and matrices are from the real field $\mathbb{R}$. $\mathbb{E}(\xi)$ denotes the expectation value of a random variable $\xi$ and $\mathbb{P}_\xi(\cdot)$ denotes the probability with respect to $\xi$. The cumulative distribution function (CDF) of $\xi$ is denoted as $\Phi_\xi$. For the optimization variables $x_t^{c,ij}$ the lower case index $t$ describes the time step, $c \in e, h$ describes the type of energy (electicity or heat) and $ij$ denote the direction of flow.

*2.2 Chance Constrained Optimization*

As in many optimization procedures, in this paper we study the following optimization problem:

$$\min_x \quad c^\top x \tag{1a}$$
$$\text{s.t.} \quad f(x,\xi) \leq 0 \tag{1b}$$
$$x \in \mathcal{X} \tag{1c}$$

where $x \in \mathbb{R}^n$ is the decision variable for which we want to find a minimum in the linear function $c^\top x$. The random vector $\xi \in \mathbb{R}^d$ only affects the constraint function $f : \mathbb{R}^n \times \mathbb{R}^d \to \mathbb{R}^m$ which in general can be vector-valued $f = (f_1, \ldots, f_m)$. The set $\mathcal{X}$ represents the feasible set of the optimization variable. As we cannot solve problem Equation 1 directly, we introduce a parameter $\epsilon$ called the violation probability.

$$\min_x \quad c^\top x \tag{2a}$$
$$\text{s.t.} \quad \mathbb{P}_\xi(f(x,\xi) \leq 0) \geq 1 - \epsilon \tag{2b}$$
$$x \in \mathcal{X} \tag{2c}$$

Now, we only guarantee compliance with the constraint $f$ with probability $1 - \epsilon = \alpha$, where $\alpha$ is called the confidence level. The formulation in Equation 2 above describes a joint probabilistic constraint. For computational reasons in many optimization frameworks, single probabilistic constraints are used by introducing $1 - \epsilon_i = \alpha_i$ for all of the $f_i$ constraint functions, thus relaxing the problem. Although the CCO problem is generally NP-hard to solve, there are some special cases which greatly simplify the problem. Firstly, in case of a single chance constraint with a linear constraint function and the uncertain variable $\xi$ following a parametric probability distribution with a "reasonably well behaved" invertible cumulative distribution function (=quantile function) $\Phi^{-1}$ (such as a multivariate gaussian), we can simplify the problem into a second order program.
Secondly, in case the uncertain variable $\xi$ does not interact directly with the optimization variable $x$, we can separate the

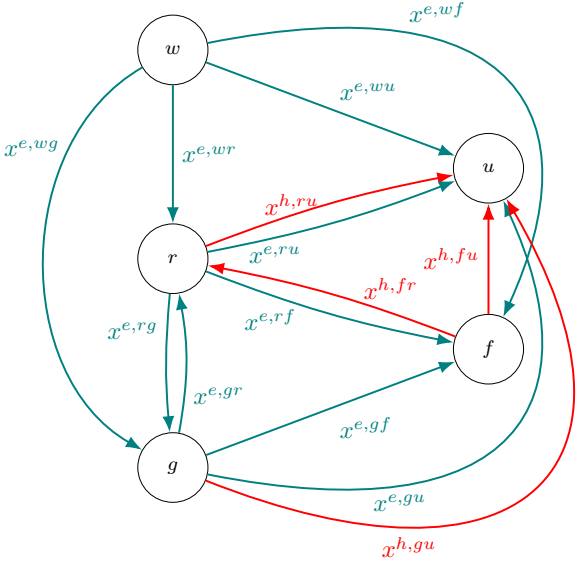

Fig. 2 Energy system model with power flow between different components. Arrows indicate variable signs. Colors indicate energy type.

deterministic and random parts of $f_i$:

$$f_i(x, \xi) = f_i^{\text{determ}}(x) + f_i^{\text{prob}}(\xi) \qquad (3)$$

Assuming further that both $f_i^{\text{determ}}(x)$ and $f_i^{\text{prob}}(\xi)$ are linear and each uncertain variable $\xi_i$ appears in exactly one constraint function $f_i$ (i.e. $d = m$), leads us to a greatly simplified linear chance constrained optimization (LCCO) problem.

$$\min_x \; c^\top x \qquad (4a)$$
$$\text{s.t.} \; a_i^\top x \leq \mu_i + \sigma_i \Phi_{\xi_i}^{-1}(1 - \epsilon_i) \qquad \forall i = 1, \ldots, m \quad (4b)$$
$$x \in \mathcal{X} \qquad (4c)$$

### 2.3 Energy System Model

We model an energy system of production, demand, storage and grid. The five positions of the system can best be understood from Figure 2. Renewable generation (mainly wind) originates in node $w$, with the quantity of generated power denoted by $E^{e,w}$. The urban and foundry nodes $u$ and $f$ have a specified demand of electric energy $D^{e,u}$ and $D^{e,u}$. The storage node $r$ can store energy either as electricity or heat with two separate storage levels $R^e$ and $R^h$. The grid can provide additional energy and heat, with prices $P^{e,g}$ and $P^{h,g}$.

The system operator's task is to manage these energy flows in a cost-effective but also reliable manner. The optimization is performed by defining values for each of the available energy flows. At every time step $t$ in the operational time horizon $T$ he sees the forecasts of the uncertain values for all time steps $t < t' \leq H$ where $H$ denotes the decision time horizon.

The optimization variables for each $t$ are

$$x_t = x_t^e + x_t^h$$
$$x_t^e = \big(x_t^{e,wu}, x_t^{e,wf}, x_t^{e,wr}, x_t^{e,wg}, x_t^{e,ru},$$
$$\qquad\qquad x_t^{e,rf}, x_t^{e,rg}, x_t^{e,gu}, x_t^{e,gf}, x_t^{e,gr}\big)$$
$$x_t^h = \big(x_t^{h,fu}, x_t^{h,gu}, x_t^{h,fr}, x_t^{h,ru}\big)$$

where

- Electric Energy Flow Variables ($x_t^e$):
  - $x_t^{e,wu}$: Electricity from *wind* to *urban* sector
  - $x_t^{e,wf}$: Electricity from *wind* to *foundry*
  - $x_t^{e,wr}$: Electricity from *wind* to *storage*
  - $x_t^{e,wg}$: Electricity from *wind* to *grid* (export)
  - $x_t^{e,ru}$: Electricity from *storage* to *urban* sector
  - $x_t^{e,rf}$: Electricity from *storage* to *foundry*
  - $x_t^{e,rg}$: Electricity from *storage* to *grid* (export)
  - $x_t^{e,gu}$: Electricity from *grid* (import) to *urban*
  - $x_t^{e,gf}$: Electricity from *grid* (import) to *foundry*
  - $x_t^{e,gr}$: Electricity from *grid* (import) to *storage*

- Heat Energy Flow Variables ($x_t^h$):
  - $x_t^{h,fu}$: Heat from *foundry* to *urban*
  - $x_t^{h,gu}$: Heat from *grid* (import) to *urban*
  - $x_t^{h,fr}$: Heat from *foundry* to *storage*
  - $x_t^{h,ru}$: Heat from *storage* to *urban*

Every decision in the decision horizon is constrained by the linear functions found in Appendix 7. These can either be enforced directly on the optimization variable $x$ or indirectly via the storage levels $R$.

The cost function can similarly be separated into parts for comprehension reasons. These functions again are entirely linear.

$$\text{TotalCost} = \sum_{t=0}^{H-1} \Big[ C_P^{e,u} \cdot \big(D_t^{e,u} - x_t^{e,wu} - \beta_e^d x_t^{e,ru} - x_t^{e,gu}\big) \quad (5a)$$
$$+ C_P^{h,u} \cdot \big(D_t^{h,u} - x_t^{h,fu} - x_t^{h,gu} - \beta_h^d x_t^{h,ru}\big) \quad (5b)$$
$$+ C_P^{e,f} \cdot \big(D_t^{e,f} - x_t^{e,wf} - \beta_e^d x_t^{e,rf} - x_t^{e,gf}\big) \quad (5c)$$
$$+ P_t^{e,g} \cdot \big(x_t^{e,gu} + x_t^{e,gf} + x_t^{e,gr}\big) \quad (5d)$$
$$+ P_t^{h,g} \cdot x_t^{h,gu} \quad (5e)$$
$$- P_t^{e,g} \cdot \big(\beta_e^d x_t^{e,rg} + x_t^{e,wg}\big) \quad (5f)$$
$$+ \text{LCOS}_e \cdot \big(x_t^{e,rf} + x_t^{e,rg} + x_t^{e,ru}\big) \quad (5g)$$
$$+ \text{LCOS}_h \cdot x_t^{h,ru} \Big] \quad (5h)$$

where

- $C_P^{e,u}$: Penalty cost per unit of unmet urban electricity demand

- $C_P^{h,u}$: Penalty cost per unit of unmet urban heat demand
- $C_P^{e,f}$: Penalty cost per unit of unmet foundry electricity demand
- LCOS$_e$: Levelized Cost of Storage for electricity — represents average cost of storing and discharging one unit of electricity
- LCOS$_h$: Levelized Cost of Storage for heat — includes capital and operational amortization per unit

The total cost function (5) accounts for all relevant economic activities across the planning horizon $t = 0, \ldots, H-1$.

The first three terms, (5a)–(5c), represent penalty costs for unmet electricity and heat demands in the urban and foundry sectors. These penalties apply to any residual demand not satisfied through local renewable production, grid purchases, or storage discharges, weighted by respective cost coefficients $C_P^{h,u}, C_P^{e,u}, C_P^{e,f}$.

Terms (5d) and (5e) capture the cost of energy purchases from the external grid. Electricity purchases include supply to the urban and foundry sectors as well as charging the electricity storage. Similarly, heat purchased from the grid supplements the urban heating demand.

Equation (5f) deducts the revenue from selling electricity back to the grid, either from discharged storage or surplus wind production. This is modeled as a negative cost, reducing the total expenditure.

The final two terms, (5g) and (5h), represent the levelized cost of storage (LCOS) for electricity and heat, respectively. These are applied proportionally to the amount of energy discharged from storage, capturing capital and operational amortization related to storage usage.

Together, the objective balances supply-demand matching, market interactions, and storage utilization in a cost-optimal manner.

*2.4 Policy formulation*

*2.4.1 Deterministic lookahead cost function approximation:*
The general idea proposed in [1] to deal with rolling forecasts is to replace an uncertain variable $\xi$ by a product of the mean $\mathbb{E}(\xi)$ and a parameter $\theta$ in order to convince the LP-algorithm to act more aggressively ($\theta > 1$) or conservatively ($\theta < 1$):

$$\text{DLA-}\theta \quad \xi \to \theta \mathbb{E}(\xi) \tag{6}$$

After setting $\theta$ we can then treat the probabilistic problem in (1) as a simple deterministic (linear) problem.

Several strategies exist to find a suitable set of parameters $\theta_{t'}$ for each time step $t'$.

*constant:* The easiest approach is to set the same $\theta_0$ for all times. This will however lead either an overestimation or an underestimation for all times.

*look-up table:* The contrary approach is to perform an optimization over all $\theta_{t'}$, however this can be numerically challenging when dealing with long time horizons $H$ and small time intervals.

*exp. decay:* As we want the LP-algorithm to act aggresssive in the far and conservative in the near future a more differentiated approach is in setting

$$\theta_{t'} = \theta^1 \exp\left(\theta^2 (t - t')\right). \tag{7}$$

Note, that the top indices to $\theta = (\theta^1, \theta^2)$ do not represent powers but are meant to distinguish from the lower case temporal index. With this formula, a combination of $\theta^1 < 1$ and $\theta^2 > 0$ will lead to the desired behaviour.

While chance-constrained optimization still transforms the problem into a linear program with a direct lookahead cost function, the formulation in (4) allows us to take into account a specific violation probability ($\epsilon$) or confidence level ($\alpha$). However, we still have the problem of defining a specific confidence level $\alpha_{t'}$ for each time step on the horizon. To accomplish this we take the same approach as in (7) by setting

$$\alpha_{t'} = \alpha^1 \exp\left(-\alpha^2 (t - t')\right). \tag{8}$$

Note the minus sign in the exponential function which leads to the confidence level decreasing in the far future, thereby enforcing a more aggressive behaviour. Once the $\alpha_{t'}$ set we can apply (4b) to receive deterministic constraint functions.

*2.4.2 Feasibility updates:* Due to the randomness of the uncertain parameters wind, demand, and price of our energy system model, it can happen that operational planning leads to infeasible situations. Errors in demand and price forecasts only leads to higher or lower values in the cost function, which is the inherent difficulty of the problem. However, being wrong about the renewable energy generation might lead to a situation where there is not enough energy available to be distributed. We therefore need to define certain rules of how to update our plan $x_{t'}^{e,w\cdot}$ before seeing the actual uncertain value $\xi = E^{e,w}$.

In the case where the realized wind generation exceeds the planned consumption, i.e.,

$$E^{e,w} > x^{e,wu} + x^{e,wr} + x^{e,wf} + x^{e,wg},$$

the system has a surplus of wind energy. Although this situation is detected in the implementation, no corrective action is currently taken.

A possible extension would be to increase the grid export $x^{e,wg}$ or fill the EES further $x^{e,wr}$ to utilize the excess generation.

Conversely, when the realized wind generation is lower than the planned consumption, i.e.,

$$E^{e,w} < x^{e,wu} + x^{e,wr} + x^{e,wf} + x^{e,wg},$$

a deficit must be resolved to maintain feasibility. This is done by reducing wind allocations in a predefined priority order: first, the grid export $x^{e,wg}$ is reduced; next, the storage charging $x^{e,wr}$ is curtailed. If further necessary, foundry demand serving $x^{e,wf}$ is reduced and compensated by increased grid consumption $x^{e,gf}$. As a last resort, urban demand serving $x^{e,wu}$ is reduced and similarly compensated by increased grid consumption $x^{e,gu}$.

## 2.5 Simulation setup

*2.5.1 Forecast Models:* In order to perform a rolling forecast optimization of the described problem we need probabilistic forecasts for the renewable energy generation, the urban electricity and heat load as well as the foundry electricity demand. We also need to provide probabilistic forecasts for the day-ahead price of electricity and heat. Although we are aware that there are highly advanced methods for these predictions, we restrict ourselves intentionally to 2-layer fully connected neural networks (FCNN) for each of these predictions. Each neural network only receives the temporal and meteorological data described in Table 1 as inputs and outputs the mean $\mu$ and standard deviation $\sigma$ of a gaussian distribution for each time steps. Time steps for one forecasts therefore don't have any inherent correlation as output by the FCNN. Instead, the correlation is enforced in a secondary step by Gaussain Process Regression.

Table 1  Neural network input features and their descriptions

| Feature | Description | Type |
| --- | --- | --- |
| year | Calendar year | $\mathbb{N}$ |
| holiday | Public holiday indicator | boolean |
| hod | Hour of day | $\mathbb{N}_0^{23}$ |
| dow | Day of week (0 = Monday) | $\mathbb{N}_0^{6}$ |
| woy | Week of year | $\mathbb{N}_0^{52}$ |
| wind_speed_bft | Wind speed (Beaufort scale) | $\mathbb{N}_0^{12}$ |
| coverage_okt | Cloud cover (in oktas) | $\mathbb{N}_0^{8}$ |
| temp_degC | Ambient temperature (°C) | $\mathbb{R}$ |

Data was gathered from several accounts. The amount of available renewable energy, the urban electrical load and the prices were forecast using NNs trained on publicly available data via the Fraunhofer ISE's `https://energy-charts.info` website. The heat demand data was taken from [19] which covers 66 households in Munich, Germany. The electrical energy demand of the foundry was provided by our project partner *Fondium AG*. Weather data was gathered from the DWD (`https://opendata.dwd.de/climate_environment/CDC/observations_germany/climate/daily/kl/historical/`) for 5 locations in Germany and subsequently averaged. During evaluation of the network the weather of only one location was used as to better reflect the actual conditions at the foundry.

As the demands and loads in the original dataset were originally collected for the entirety of Germany we adjusted the data by the number of households in Germany and the number of households in the urban area of the energy system

*2.5.2 Implementation Details:* For this case study we studied the cost accumulated throughout one year (2023) with hourly rolling with a forecast and optimization horizon $H = 24$. We used 10 samples for every configuration of the hyperparameters $\theta$ and $\alpha$ and for the selection of appropriate storage sizes. Details on the parameters of the energy system model can be found in Table 2 The python PuLP package (`https://coin-or.github.io/pulp/#`) was used for the exact optimization process.

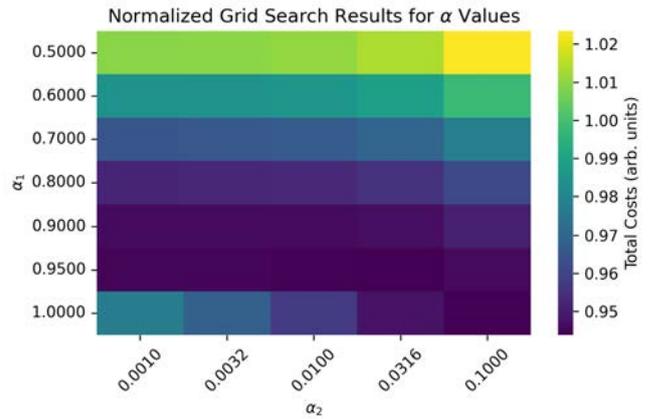

Fig. 3  Normalized costs from the grid search over the Parameter set $\alpha = (\alpha^1, \alpha^2)$ for the base scenario

## 3 Results

We evaluate the policies introduced in subsection 2.4 on the case study introduced in subsection 2.5. Specifically, we show grid search results for the best set of parameters for determining forecast values to be used by the lookup table and we investigate the ideal size of storages to be used.

### 3.1 Policy Parameter Search

As can be seen in Figure 4 and Figure 3 a grid search was performed to define the best parameter sets $\theta^*$ and $\alpha^*$. We plot the evolution of $\theta_{t'}$ and $\alpha_{t'}$ in Figure 5. It clearly shows for both that a conservative approach in planning for energy dispatch is favorable in general.

For $\theta$ the exact value $\theta^1$ seems to be of minor importance leading to our belief, that this simple model cannot deal with uncertainties in the far future and thus only acts in short term planning.

For $\alpha$ we can see that the model exactly captures our desired behaviour, using quite conservative confidence levels ($> 0.9$) in the short range and quite aggressive confidence levels ($< 0.5$) in the long range.

### 3.2 Storage Size Search

The grid search over possible storage sizes was performed for the base case with no adaptation of the forecasts. The results can be seen in Figure 6.

In can clearly be seen that the relative size of the EES ($sf_e$) is much more important than the relative size of the TES ($sf_h$) in determining the overall costs. We also expected the optimal EES and TES size to correlate with one another; however, the opposite seems to be the case. They are very much independent of each other.

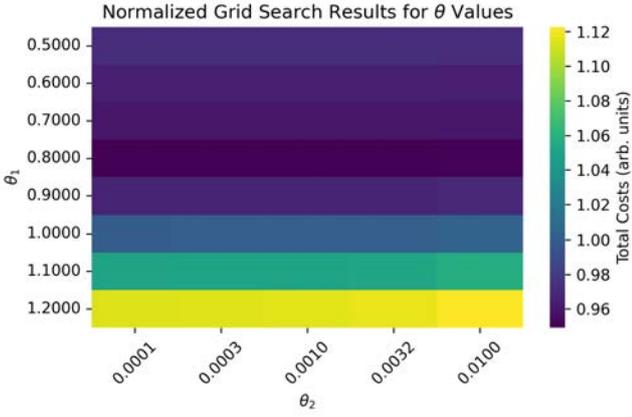

Fig. 4 Normalized costs from the grid search over the Parameter set $\theta = (\theta^1, \theta^2)$ for the base scenario

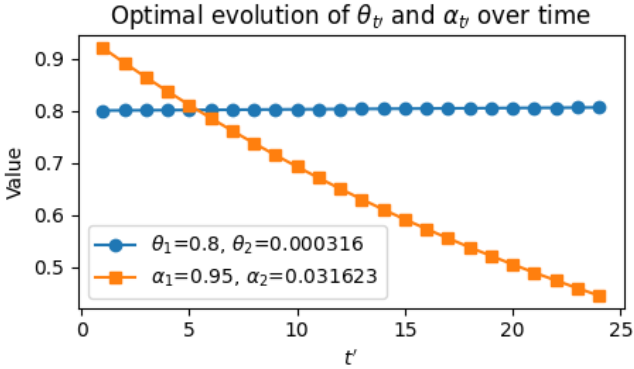

Fig. 5 Optimal evolution of $\theta_{t'}$ (blue) and $\alpha_{t'}$ (orange) over the horizon according to the results obtained in the parameter grid searches.

Still, the optimal storage factors as determined by this grid search are $sf^* = (sf_e, sf_h) = (10^{-3}, 10^0)$, which corresponds to an actual capacity of only $R_e^{\max*}, R_h^{\max*} = (100\,\text{kWh}, 100\,000\,\text{kWh})$. We have the following thoughts on these results: The relatively low value of $R_e^{\max*}$ stems from the fact that the levelised cost of storage has been modelled using a number of approximate values.. Also, the excess wind energy could not be used in our model and therefore further reduces incentives for a large EES. Still, these results show that the usage of a large production plant with heat usage can be highly beneficial.

The reason for the low influence of storage size on cost might be the steady price of grid heat, which was assumed. Even though this is realistic for gas prices on a long time scale, as urban heating turns to central heating with heat pumps or hydrogen, we can assume these prices also become more volatile.

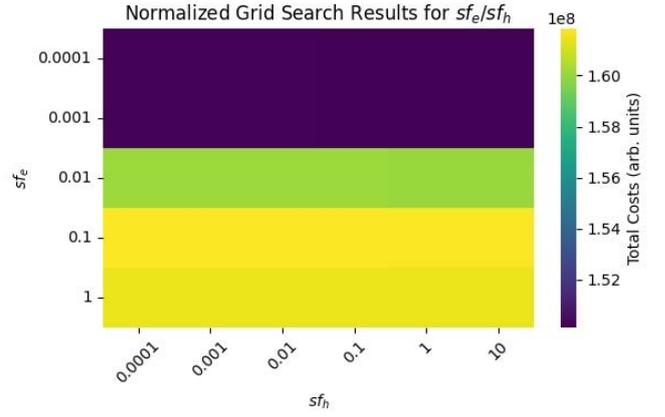

Fig. 6 Network operation cost associated with relative storage sizes of EES (y-axis) and TES (x-axis). The storage factors $sf$ refer to the capacity of EEs and TES specified in table Table 2

## 4 Conclusion

This paper establishes a risk-aware strategy for dealing with rolling forecasts in an energy network with a foundry as a heat provider for a distributed energy network. It estimates the yearly benefit of using these rolling forecasts with an optimized confidence level and gives a recommendation on storage size.

It shows once again the importance of including advanced methods for forecasting electrical and thermal loads and generation mechanisms as high confidence levels are important in order to make maximum profit.

Future avenues for research include more advanced models for the uncertainty estimation as well as the inclusion of uncertainty in more parameters of the energy system. We emphasize the need for models that especially capture the tail behaviour of the probability distributions well, as these are important for high confidence levels. The models should also incorporate some autoregressive components in the future. We also recommend the investigation of a method for controlling the foundry output and of the transferability of these results from the financial to the environmental realm by including CO2 costs in the yearly costs as well.

## 5 Acknowledgment

The author(s) gratefully acknowledge the support of the Carl-Zeiss Stiftung, whose funding of the project "DeepCarbPlanner" made this research possible.

## 7 Appendix: Constraint Functions of the Energy System Model

$$x_t^{e,wu} + x_t^{e,wf} + x_t^{e,wr} + x_t^{e,wg} \leq E_t^{e,w} \tag{9a}$$

$$x_t^{e,ru} + x_t^{e,rf} + x_t^{e,rg} \leq \min(\gamma_e^d, R_t^e) \tag{9b}$$

$$\beta_e^c(x_t^{e,wr} + x_t^{e,gr}) \leq \gamma_e^c \tag{9c}$$

$$x_t^{e,wu} + \beta_e^d x_t^{e,ru} + x_t^{e,gu} \leq D_t^{e,u} \tag{9d}$$

$$x_t^{e,wf} + \beta_e^d x_t^{e,rf} + x_t^{e,gf} \leq D_t^{e,f} \tag{9e}$$

$$x_t^{h,ru} \leq \min(\gamma_h^d, R_t^h) \tag{9f}$$

$$\beta_h^c x_t^{h,fr} \leq \gamma_h^c \tag{9g}$$

$$x_t^{h,fu} + x_t^{h,gu} + \beta_h^d x_t^{h,fu} \leq D_t^{h,u} \tag{9h}$$

$$x_t^{h,fu} + x_t^{h,fr} \leq \delta^{eh}(x_t^{e,wf} + \beta_e^d x_t^{e,rf} + x_t^{e,gf}) \tag{9i}$$

$$R_{t+1}^e = R_t^e + \beta_e^c(x_t^{e,wr} + x_t^{e,gr}) - (x_t^{e,ru} + x_t^{e,rf} + x_t^{e,rg}) \tag{10a}$$

$$R_{t+1}^h = R_t^h + \beta_h^c x_t^{h,fr} - x_t^{h,ru} \tag{10b}$$

$$R_{t+1}^e \leq R_e^{\max} \tag{10c}$$

$$R_{t+1}^h \leq R_h^{\max} \tag{10d}$$

where several constant parameters appear

- $\beta_e^c$: Charging efficiency of the EES
- $\beta_e^d$: Discharging efficiency of the EES
- $\beta_h^c$: Charging efficiency of the TES
- $\beta_h^d$: Discharging efficiency of the TES
- $\gamma_e^c$: Maximum electric charging rate
- $\gamma_e^d$: Maximum electric discharging rate
- $\gamma_h^c$: Maximum heat charging rate
- $\gamma_h^d$: Maximum heat discharging rate
- $\delta^{eh}$: Efficiency of usable heat conversion at the foundry

The constraints described in equations (9a)–(10d) correspond to the fundamental inventory dynamics of an integrated energy system. Equation (9a) ensures that the total electricity drawn from wind—whether it is directed toward urban and foundry demand, storage charging, or grid export—does not exceed the available wind energy supply. The discharge processes for both the electrical energy storage (EES) and thermal energy storage (TES) are governed by constraints (9b) and (9f), which cap the discharge based on each system's current storage level and discharge efficiency. Similarly, constraints (9c) and (9g) restrict the charging rates of EES and TES, incorporating both the system's technical charging limits and associated efficiencies. The system's ability to meet energy demand is addressed in constraints (9d) and (9e), which ensure that electricity needs in the urban and foundry sectors are fulfilled by a combination of wind generation, storage discharge, and grid purchases. For heat provision, constraint (9h) guarantees that urban heat demand is similarly satisfied through foundry waste heat, grid imports, and discharged thermal storage. The coupling constraint (9i) limits the total amount of heat that can be produced by the foundry, ensuring it does not exceed what can be thermodynamically supported by the electricity used for its operation, scaled by a conversion efficiency.

Storage dynamics are explicitly modeled in equations (10a) and (10b), capturing the evolution of the storage levels from one time step to the next based on the net result of charging and discharging. Finally, constraints (10c) and (10d) impose upper bounds on the storage capacities of the EES and TES, ensuring that the stored energy never exceeds their physical limits.

## 8 Appendix: Detailed Parameters of the Energy System Model

Table 2 Energy Model parameters

| Node | Symbol | Description | Value | Unit |
|---|---|---|---|---|
| urban | $n_{\text{houses}}$ | number of households | 10000 | - |
| EES | $R_e^{\max}$ | Max storage capacity | $10 \cdot n_{\text{houses}}$ | kWh |
|  | $R_e^{\min}$ | Min storage level | 0 | kWh |
|  | $R_e^0$ | Initial storage level | 0 | kWh |
|  | $\beta_e^c$ | Charging efficiency | 0.9 | – |
|  | $\beta_e^d$ | Discharging efficiency | 0.9 | – |
|  | $\gamma_e^c$ | Max charge rate | $R_e^{\max}/10$ | kW |
|  | $\gamma_e^d$ | Max discharge rate | $R_e^{\max}/10$ | kW |
|  | $\text{LCOS}_e$ | Levelized cost of storage | 0.05 | €/kWh |
| TES | $R_h^{\max}$ | Max storage capacity | $10 \cdot n_{\text{houses}}$ | kWh |
|  | $R_h^{\min}$ | Min storage level | 0 | kWh |
|  | $R_h^0$ | Initial storage level | 0 | kWh |
|  | $\beta_h^c$ | Charging efficiency | 0.9 | – |
|  | $\beta_h^d$ | Discharging efficiency | 0.9 | – |
|  | $\gamma_h^c$ | Max charge rate | $R_h^{\max}/10$ | kW |
|  | $\gamma_h^d$ | Max discharge rate | $R_h^{\max}/10$ | kW |
|  | $\text{LCOS}_h$ | Levelized cost of storage | 0.01 | €/kWh |
| Foundry | $\delta^{eh}$ | Electr.-to-heat efficiency | 0.5 | – |